\documentclass[10pt]{amsart}
\usepackage{latexsym, amsmath, amssymb,amsthm,amsopn,amsfonts,amscd}
\usepackage{version}
\usepackage{epsfig,graphics,color,graphicx,graphpap}

\newcommand{\RR}{{\mathbb R}}

\newcommand{\nn}{\nonumber}

\theoremstyle{plain}

\newtheorem{thm}{Theorem}
\newtheorem{prop}{Proposition}[section]

\newtheorem{lem}[prop]{Lemma}

\theoremstyle{definition}

\newtheorem{rem}{Remark}[section]
\newtheorem{defn}[prop]{Definition} 

\numberwithin{equation}{section}

\def\squarebox#1{\hbox to #1{\hfill\vbox to #1{\vfill}}}

\newcommand{\p}{\partial}
\usepackage{amsxtra}

\usepackage{fancyhdr}
\ifx\pdfoutput\undefined
  \DeclareGraphicsExtensions{.eps}
\else
  \ifx\pdfoutput\relax
    \DeclareGraphicsExtensions{.eps}
  \else
    \ifnum\pdfoutput>0
      \DeclareGraphicsExtensions{.pdf}
    \else
      \DeclareGraphicsExtensions{.eps}
    \fi
  \fi
\fi

\title
[Wave operator bounds]
{Wave operator bounds for $1$-dimensional Schr\"odinger operators with
singular potentials and applications}  

\author[V. Duch\^ene]
{Vincent Duch\^ene}
\email{vincent.duchene@ens.fr}

\address{\'Equipe EDP, DMA - \'Ecole Normale Supe\'rieure \\
45, rue d'Ulm, 75230 Paris Cedex 05 - France }

\author[J.L. Marzuola]
{Jeremy L. Marzuola}
\email{jm3058@columbia.edu}
\address{Department of Applied Physics and Applied Mathematics, Columbia University \\
200 S. W. Mudd, 500 W. 120th St., New York City, NY 10027, USA}
\author[M.I. Weinstein]
{Michael I. Weinstein}
\email{miw2103@columbia.edu}
\address{Department of Applied Physics and Applied Mathematics, Columbia University \\
200 S. W. Mudd, 500 W. 120th St., New York City, NY 10027, USA}

\begin{document}    

\begin{abstract} 
Boundedness of wave operators for Schr\"odinger operators in one space dimension for 
a class of singular potentials, admitting finitely many Dirac delta distributions,  is proved. Applications are presented to, for example, dispersive estimates and commutator bounds.
 \end{abstract}

\maketitle

%\tableofcontents

\section{Introduction}
\label{sec:intro}

{\it Wave operators} provide a means for converting operator bounds for a ``free'' dynamics generated
 by a constant coefficient Hamiltonian, $H_0=-\Delta$ to analogous operator bounds about ``interacting"
  dynamics associated with a variable coefficient Hamiltonian, $H=-\Delta+V$, on its continuous spectral subspace. Indeed let $W_{\pm}$ and $W_{\pm}^*$ denote  wave operators associated with the free and interacting Hamiltonians $H_0$ and $H$
   (defined by \eqref{eqn:wpm} and \eqref{eqn:wpmstar}). Then we have
\begin{align}
W_{\pm} W_{\pm}^* &= P_c, \ \ \ W_{\pm}^* W_{\pm} = Id \label{ww-star}\\
f(H) P_c &= W_{\pm} f(H_0) W_{\pm}^*, \ f(H_0) = W_{\pm}^* f(H)
W_{\pm},\ \ \ f\ {\rm Borel\ on}\ \RR\ .\label{wfwstar}
\end{align}
It follows that bounds on $f(H)P_c$ acting between $W^{k_1,p_1}(\RR^d)$ and $W^{k_2,p_2}(\RR^d)$ can be derived from bounds on $f(H_0)$ between these spaces if the wave operators $W_{\pm}$ are bounded between  $W^{k_1,p_1}(\RR^d)$ and $ W^{k_2,p_2}(\RR^d)$ for $k_j\ge0$ and $p\ge1$. Here, $W^{k,p}(\RR^d),\ k\ge1,\ p\ge1$ denotes the Sobolev space of functions having derivatives up to order $k$ in $L^p(\RR^d)$. 

Applications along the lines of the above discussion have appeared in \cite{MW}. For example, 
\begin{equation}
\left\| e^{-iHt}P_c(H) f \right\|_{L^p(\RR^d)}\ =\ \left\| W_{\pm} e^{-iH_0t} W_{\pm}^* f\right\|_{L^p(\RR^d)}\ \le\ C\ |t|^{-\frac{d}{2}-\frac{d}{p}}\ \left\| f\right\|_{L^q(\RR^d)},\ \ \ p^{-1}+q^{-1}=1,\ \ p\ge1.
\nn\end{equation}

Boundedness of wave operators in $W^{k,p} (\RR^d)$, under smoothness and decay assumptions on $V(x)$ was proved in \cite{Ya} in dimensions $d\ge2$.  Weder \cite{W} proved boundedness in dimension one; see also \cite{DF}. 
In \cite{W} it is assumed that $V\in L^1_\gamma (\RR)$, the space of all complex-valued measurable functions $\phi$ defined on $\RR$ such that
\begin{eqnarray}
\| \phi \|_{L^1_\gamma} = \int | \phi (x) | (1 + |x|)^\gamma dx <
\infty .
\end{eqnarray}
For $V$ falling into a class of generic potentials, the assumption is $\gamma > 3/2$, otherwise it is assumed $\gamma > 5/2$.
 
 Schr\"odinger operators with singular potentials arise in models, which have recently been extensively investigated. See, for example,
\cite{MW,GHW,JW,HMZ,HMZ1,HZ}, where Dirac delta function potentials are considered.  Boundedness of wave operators for singular potentials satisfying the hypotheses of Theorem \ref{thm:main} is used implicitly in references \cite{JW} and \cite{GHW}, but this boundedness appears not to have been addressed previously.  This gap in the literature is addressed by the current work.  Another motivation for the present work is the study of scattering for highly oscillatory structures in the homogenization limit \cite{DucheneWeinstein}, where  bounds on  $(m^2+H)^{-1}P_c(H)(m^2-\partial_x^2)$, where $H= -\partial_x^2+V(x)$ is a Schr\"odinger operator with a singular (distribution) part to the potential $V(x)$, are required; see section \ref{sec:nonapp}.

  This article is devoted to an extension of the one-dimensional results \cite{W} to the case of singular potentials. In particular, our results apply to Hamiltonians of the form 
\begin{equation} H\ =\ -\partial_x^2\ +\ V(x), \nn\end{equation}
where $V(x)$ satisfies: 
\bigskip

\noindent{\bf Hypotheses (V)}
\begin{align}
\label{Vdecomp}
V(x)\ &=\ V_{sing}(x)\ +\ V_{reg}(x),\\
V_{sing}(x)\ &=\  \sum_{j=0}^{N-1}\ q_j\ \delta (x-y_j),\ \ 
 q_j, y_j \in \RR,\ \ y_j < y_{j+1},\ \ q_j \neq 0\label{Vsing},\\
\|V\|_{L^1_{{3\over2}+}(\RR)}\ &\equiv\ \int_\RR (1+|s|)^{{3\over2}+} |V_{reg}(s)|\ ds\ <\ \infty.
\label{Vreg}\end{align}
\bigskip

The paper is structured as follows. In section \ref{sec:main} we state our main result, Theorem \ref{thm:gen}, concerning boundedness of wave operators. In section \ref{sec:strategy} the strategy of proof is outlined. Section \ref{sec:back} summarizes facts
 about Jost solutions, distorted plane waves, reflection and transmission coefficients {\it etc}.   Some related technical results  are contained in Appendix \ref{app:mjbounds}. In section \ref{sec:state} we state a general result, Theorem \ref{thm:main}, from which Theorem \ref{thm:gen} follows. The proof of Theorem \ref{thm:main} is given in section \ref{sec:proof-of-Whighlow}. Finally, in section \ref{sec:nonapp} we present examples (multi- delta function potentials) and applications to dispersive estimates, commutator bounds and well posedness. 
 \bigskip
 
{\bf Acknowledgements:}\ 
JLM was supported, in part, by a  U.S. National Science Foundation Postdoctoral Fellowship   in the Department of Applied Physics and Applied Mathematics (APAM) at Columbia University.  MIW  was supported, in part, by  U.S. NSF Grant DMS-07-07850.
JLM and MIW wish to acknowledge the hospitality of the Courant Institute of Mathematical Sciences, where MIW was on sabbatical during the preparation of this manuscript. VD was supported, in part, by Agence Nationale de la Recherche Grant ANR-08-BLAN-0301-01. VD would like to thank APAM for its hospitality during the Spring of 2008, when this work was initiated. 
 
\section{Main results}\label{sec:main}

We first define and review properties of the wave operators.
For basic results on wave operators see, for  example, \cite{Ag,RSv3,Schec}.

Introduce the self-adjoint operators $H_0 = -\Delta$ and  $H = -\Delta
+ V$. Here, $V$ is a real-valued potential, satisfying assumptions
given below; see Section \ref{sec:state}.  Let $P_c=P_c(H)$ denote the continuous spectral projection associated with $H$.
The wave operators, $W_\pm$ and their adjoints $W_\pm^*$ are defined by
\begin{align}
\label{eqn:wpm}
W_\pm &\equiv s - \lim_{t \to \infty} e^{it H} e^{-it H_0}\\
W_\pm^* &\equiv s - \lim_{t \to \infty} e^{i t H_0} e^{-it H} P_c.
\label{eqn:wpmstar}
\end{align}
The wave operators satisfy the properties \eqref{ww-star} and \eqref{wfwstar}.
The notion of wave operators is intimately related to the idea of
distorted Fourier bases, which are discussed in detail in \cite{Ag},
\cite{Ho2}, \cite{RSv4}.  In one dimension, this is directly related
to the Jost solutions.  These objects are studied in general in
\cite{RSv4} and generalized to even a certain class of
non-self-adjoint operators in \cite{KS}.

Our main result, Theorem \ref{thm:main}, combined with the calculations of Section \ref{sec:singpot}, implies the following:
\begin{thm}
\label{thm:gen}
Consider the Schr\"odinger operator with a potential, $V(x)$, satisfying {\bf  Hypotheses (V)}.
 Then $W_\pm$ and
$W^*_\pm$ originally defined on $W^{1,p} \cap L^2$, $1 \leq p \leq
\infty$, have extensions to bounded operators on $W^{1,p}$, $1 < p <
\infty$.  Moreover, there are constants $C_p$ such that:
\begin{eqnarray}
\| W_\pm f \|_{W^{1,p}(\RR)} \leq C_p \| f \|_{W^{1,p}(\RR)}, \ \| W^*_\pm f \|_{W^{1,p}(\RR)} \leq C_p \| f \|_{W^{1,p}(\RR)}, \ f  \in W^{1,p}(\RR), \ 1 < p < \infty.
\end{eqnarray}
\end{thm}

\begin{rem} In general, the wave operators are not bounded on $L^1$. The constraint  $p>1$ is due to the Hilbert transform, $\mathcal{H}$ not being bounded on $L^1$; see \cite{W}.
\end{rem}

\section{Strategy of Proof}
\label{sec:strategy}

We use the approach for wave operators on $\RR$ initiated by Weder in \cite{W}. 
 The heart of the matter concerns the detailed low and high frequency behavior of Jost solutions, worked out by Deift and Trubowitz \cite{DT}, or a consequence of their methods. 
  The idea  is to split the wave operators into high and low
frequency components:
\begin{eqnarray*}
W_{\pm} = W_{\pm,high} + W_{\pm,low}.
\end{eqnarray*}

For the high frequency component we prove for $\phi \in \mathcal{S}$,
\begin{eqnarray*}
W_{\pm,high} \phi = \sum_j S_{A_j} \phi,\ \ {\rm where}\ \ 
S_A \phi \equiv \int_{-\infty}^\infty A(x,y) \phi (y) dy.
\end{eqnarray*}
For each $A=A_j$, we use the criterion  (Young's inequality \cite{Folland}) for $L^p,\ 1\le p\le\infty$ boundedness:
\begin{align}
C_A\ &\equiv\ \sup_{x\in\RR}\ \int_\RR\ |A(x,y)|\ dy\ +\ \sup_{y\in\RR}\ \int_\RR |A(x,y)|\ dx\ <\ \infty\nn\\
&\implies \left\|S_A\phi\right\|_{L^p}\ \le\ C_A\ \left\|\phi\right\|_{L^p}.\nn
\label{young}\end{align}
to prove
\begin{equation}
\left\| W_{\pm,high} \phi \right\|_{W^{k,p}}\ \le\ C_p\ \left\|  \phi \right\|_{W^{k,p}},\ \ 1<p<\infty,\ \ k\ge0.
\label{Whigh}
\end{equation}

For the low frequency components, we have
\begin{eqnarray*}
W_{\pm,low} \sim \mathcal{H}\ + \  \sum_j S_{A_j},
\end{eqnarray*} 
where $S_{A_j}$ is as above and 
 $\mathcal{H}$ denotes the Hilbert Transform 
\begin{eqnarray}\label{HTdef}
(\mathcal{H} \phi) (x) & = & \frac{1}{\pi}\ {\rm P.V.}\ \int \frac{\phi(x-y)}{y} dy\ =\   \int_{-\infty}^\infty e^{ikx}\ \left(-i\ {\rm sgn}(k)\right)\ \hat{\phi} (k) dk
\end{eqnarray}
Here,  $F$ and $F^{-1}$ denote the Fourier Transform on $\RR$ and its inverse,  defined by
\begin{equation}
\hat{\phi}(k)\ \equiv\ F \phi(k) = \frac{1}{2 \pi} \int e^{-ikx} \phi (x) dx,\ \ \ 
\check{\Phi}(x)\ \equiv\ F^{-1}\Phi(x) = \int e^{ikx} \Phi (k) dk.
\end{equation}
 Thus, for low frequencies, boundedness 
 \begin{equation}
\left\| W_{\pm,low} \phi \right\|_{W^{k,p}}\ \le\ C_p\ \left\|  \phi \right\|_{W^{k,p}},\ \ 1<p<\infty,\ \ k\ge0
\label{Wlow}
\end{equation}
 reduces to the boundness properties of the Hilbert  
  transform \cite{Stein}:
\begin{thm}\label{thm:CZest}
$\mathcal{H} : W^{k,p} \to W^{k,p}, {\rm for}\ 1<p<\infty\ \ {\rm  and} \ \ k\ge0,$
with $\left\| \mathcal{H}\phi \right\|_{W^{k,p}(\RR)}\ \le\ K_p\ \| \phi \|_{W^{k,p}(\RR)}$.
\end{thm}
Estimates \eqref{Whigh} and \eqref{Wlow} then imply the theorem.  The proof of  \eqref{Whigh} and \eqref{Wlow} is given in section \ref{sec:proof-of-Whighlow}.
 We now develop some background for  implementing the strategy. 

\section{Background spectral theory of $H=-\partial_x^2+V$}
\label{sec:back}

\subsection{Distorted plane waves, $e_\pm(x;k)$}
Consider the operator $H = - \p_x^2 + V(x)$, defined as a self-adjoint
operator on $L^2 (\RR)$.
   Denote by $P_d$ and $P_c$ the discrete and continuous spectrum projections. $P_d$ and $P_c$ are orthogonal projections with $P_c = Id - P_d$.
   
   Denote by $R_0$  the outgoing ``free'' resolvent operator $R_0(k)=(-\p_x^2 - k^2)^{-1}$
 with kernel 
 \begin{equation}
 R_0(k)(x,y)=-(2ik)^{-1}\exp(ik|x-y|)
 \nonumber\end{equation}
 and finally introduce the {\it distorted plane waves}, $e_\pm(x;k)$:

\begin{defn}
$u=e_\pm (x;k)$ are the unique solutions to $(H-k^2) u =0$ satisfying 
\begin{equation}
e_\pm(x;k)\ =\  e^{\pm i k x}\ +\ {\rm outgoing}(x),
\label{outgoing}
\end{equation}
where a function 
$U$ is said to be outgoing as $|x|\to\infty$ if
\begin{equation}
\left(\ \partial_x\ \mp\ ik\ \right)U =\ 0,\  x\to\pm\infty.
\nn\end{equation}
Thus, $e_\pm(x,;k)$ is given by the integral equation:
\begin{eqnarray*}
e_{\pm} (x;k) = e^{\pm i k x} - R_0(k) V e^{\pm i k x}.
\end{eqnarray*}
\end{defn}

The continuous spectral projection, $P_c$, is given by
\begin{eqnarray}
P_c f (x) =  \frac{1}{2\pi} \int \int_0^\infty \left( e_+(x,k)\ \overline{e_+ (y,k)} + 
e_- (x,k)\ \overline{e_- (y,k)} \right) f(y) dk dy.
\label{Pcf}\end{eqnarray}
see, for example, \cite{TZ}.

We write
\begin{align}
P_cf\ &\equiv \ F_+^*\ F_+\ f,\ \ {\rm where\ it\ follows\ from\ \eqref{Pcf}\ that }\ \nonumber\\
F_+f\ &\equiv\ \int_\RR\ \overline{\Psi_+(y,k)}\ f(y)\ dy,\ \ \ \ F^*_+f\ \equiv\ \int_\RR\ {\Psi_+(y,k)}\ f(y)\ dy\ \ {\rm and}\ \ \label{F+def}\\
\Psi_+ (y,k) & =  \frac{1}{\sqrt{2 \pi}} 
\left\{ \begin{array}{cc}
{e}_{+} (x;k) & k \geq 0, \\
{e}_{-} (x;-k) & k < 0 
\end{array} \right. \label{Psi+def}
\end{align}
We also define  $\Psi_{-} (x,k) = \overline{\Psi_+ (x,-k)}$.  

\subsection{Jost solutions}
To make direct use of the arguments in \cite{W} and \cite{DT}, we express the results of the preceding subsection in terms of {\it Jost solutions}, commonly introduced for one-dimensional Schr\"odinger operators.

Given the Schr\"odinger equation
\begin{eqnarray}
\label{eqn:statsch}
- \frac{d^2}{dx^2} u + V u = k^2 u, \ k \in \mathbb{C},
\end{eqnarray}
we define the Jost solutions, $f_j(x,k)$, $j = 1,2$, $\text{Im} k \geq 0$,
 to be the  unique solutions of \eqref{eqn:statsch} satisfying the conditions:
 \begin{align}
 f_1 (x,k) &\ -\  e^{ikx}\to0,\ \ \ x \to \infty,\ \ {\rm and}\nonumber \\
f_2 (x,k) &\ -\ e^{-ikx}\to0,\ \ \ x \to -\infty. \label{Jostdef}
\end{align}
The Jost solutions are linearly  independent solutions of 
\eqref{eqn:statsch} for $k \neq 0$. Therefore, there are unique functions
$T(k)$, $R_j (k)$, $j=1,2$ such that for $k \in \mathbb{R} \setminus 0$
\begin{eqnarray}
\label{eqn:f1f2}
f_2 (x,k) & = & \frac{R_1(k)}{T(k)}\ f_1 (x,k) + \frac{1}{T(k)}\ f_1 (x,-k), \\
f_1 (x,k) & = & \frac{R_2(k)}{T(k)}\ f_2 (x,k) + \frac{1}{T(k)}\ f_2 (x,-k)
\end{eqnarray}
For a potential, $V$, with compact support within $(-r,r)$,  $R_j(k)$ and $T(k)$ are defined via the solutions:
\begin{align}
u_1(x;k) &= \left\{ \begin{array}{cc}
e^{ikx}+R_2(k)e^{-ikx}, & x<-r, \label{u1exp} \\ 
T(k)e^{ikx}, &x>r\\
\end{array}\right.\\
u_2(x;k) &= \left\{ \begin{array}{cc}
T(k)e^{-ikx}, & x<-r,\\ 
e^{-ikx}+R_1(k)e^{ikx}, &x>r \label{u2exp} \\
\end{array}\right.
\end{align}
Generically, 
\begin{equation}
T(k) = \alpha k
+ o(k),\ \ \ \ 1+R_j (k) = \alpha_j k + o(k),\ \ \ \ j=1,2,\ \ \ k \to 0.
\end{equation}  
$T(k)$ is called the transmission coefficient associated with $H$. $R_1 (k)$ is the right to left reflection coefficient, and $R_2(k)$
the left to right reflection coefficient.  

It follows from  \eqref{outgoing}, \eqref{Jostdef} and \eqref{eqn:f1f2} that 
\begin{eqnarray}
\Psi_+ (x,k) 
& = &  \frac{1}{\sqrt{2 \pi}} \left\{ \begin{array}{cc}
T(k)\ e^{ikx}\ m_1 (x,k) & k \geq 0, \\
T(-k)\ e^{ikx}\ m_2 (x,-k) & k < 0,
\end{array} \right. 
\label{Psi+T}\end{eqnarray}
where $m_1(x,k)-1 \to\ 0$ as $x\to\infty$ and $m_2(x,k)-1\ \to\ 0$ as $x\to-\infty$.
 The detailed smoothness and decay properties, in $x$ and $k$, of $m_j(x;k)-1$ are 
  required in estimates. These are given in Appendix \ref{app:mjbounds}.

\section{Statement of the Main Theorem}
\label{sec:state}

Our main result, from which Theorem \ref{thm:gen} follows, is:  
\begin{thm}
\label{thm:main}
Let  $H = -\p_x^2 + V(x)$ be self-adjoint on $L^2
(\RR)$, where $V=V_{sing}(x)+V_{reg}(x)$ for which   the transmission and reflection coefficients (see \eqref{eqn:f1f2}) satisfy the bounds:
\begin{eqnarray}
|R(k)|, \ |T(k)-1|, \ |\partial_kR(k)|, \ |\partial_k T(k)| \leq \frac{C}{\langle k
  \rangle}. 
  \label{RT-assume}
  \end{eqnarray}
  Assume further that there exists $a>0$ sufficiently large such that
   \begin{eqnarray}
   | \p_x^\alpha m_1 (x,k)| \!\!\!\!\! &+& \!\!\!\!\! | \p_x^\alpha m_2 (x,k)| \leq C(a) \ \text{for} \ |x| \leq a, \ \alpha = 0,1,  \label{m1m2compbound} \\
\left|\ m_1(x;k)-1\  \right|\ \!\!\!\!\! &+& \!\!\!\!\! \  \left|\ \partial_k m_1(x;k)\ \right|\ +\
 \left|\ \partial_x m_1(x;k)\ \right| \le\ C\ \frac{\int_x^\infty |V_{reg}(t)| (1+|t|)\ dt}{1+|k|},\ \ \ x\ge a, \label{m1bound}\\
 \left|\ m_2(x;k)-1\  \right|\ \!\!\!\!\! &+& \!\!\!\!\! \  \left|\ \partial_k m_2(x;k)\ \right|\ +\
 \left|\ \partial_x m_2(x;k)\ \right| \le\ C\ \frac{\int_{-\infty}^x |V_{reg}(t)| (1+|t|)\ dt}{1+|k|},\ \ \ x\le -a\ \ .\label{m2bound}
\end{eqnarray}
Then $W_\pm$ and $W^*_\pm$ originally defined on $W^{1,p} \cap L^2$, $1 \leq p \leq \infty$,  extend to bounded operators on $W^{1,p}$, $1 < p < \infty$.  Furthermore, there are constants $C_p$ such that:
\begin{eqnarray}
\| W_\pm f \|_{W^{1,p}} \leq C_p \| f \|_{W^{1,p}}, \ \| W^*_\pm f \|_{W^{1,p}} \leq C_p \| f \|_{W^{1,p}}, \ f  \in W^{1,p} \cap L^2, \ 1 < p < \infty.
\end{eqnarray}
\end{thm}

\begin{rem}
Deift and Trubowitz \cite{DT} establish the bounds  \eqref{m1bound} and \eqref{m2bound}   for any potential $V(x)$, for which $(1+|x|)\left|V(x)\right|\in L^1(\RR)$ with $a=0$. Their proof applies to a potential of the type in {\bf Hypothesis (V)}, $V=V_{sing}+V_{reg}$, where  $V_{sing}$ has a finite set of Dirac masses within  an interval $(-A,A)$, and such that  $(1+|x|)\left|V_{reg}(x)\right|\in L^{{3\over2}+}(\RR)$. In this case the bounds \eqref{m1bound} and \eqref{m2bound} hold with $a=A$,  $C$ depending on $A$ and $V$ replaced by $V_{reg}$.
\end{rem}

\begin{rem}
In fact, less restrictive bounds on $V_{reg}$ as developed in \cite{DF} would suffice.  However, for simplicity we will follow the work of \cite{W} as it makes some computations more explicit.
\end{rem}

\section{Proof of Main Theorem \ref{thm:main}}
\label{sec:proof-of-Whighlow}

We follow the strategy described in section \ref{sec:strategy}. 

Let $\chi(x\ge1) \in C^\infty (\RR)$ denote non-decreasing cut-off functions
such that 
\begin{eqnarray}
\chi(x\ge1) =
 \left\{ \begin{array}{cc}
   0\ & x \le \frac{1}{2}, \\
 1\   & x \ge 1\ .
\end{array} \right. 
\end{eqnarray}

To localize in frequency space, introduce $\  \psi\left(|k|\le k_0\right) \in C_0^\infty (\RR)$ be a compactly
supported cut-off function, depending on a parameter, $k_0$, to be chosen,
 such that 
\begin{eqnarray}
\psi\left(|k|\le k_0\right) =
 \left\{ \begin{array}{cc}
   1\ & |k| \le  \ k_0, \\
 0\   & |k| \ge\ 2k_0 \ \ .
\end{array} \right. 
\end{eqnarray}

We decompose any $\phi\in L^2(\RR)$ into its low and high frequency parts:
\begin{equation}
\phi(x)\ =\ \phi_{low}(x)\ +\ \phi_{high}(x),\ \ \ \ {\rm where\ using\  }\ \ D\equiv-i\partial_x,
\label{low+high}\end{equation}
\begin{align}
\phi_{low}(x)\ &\equiv\ \psi(|D|\le k_0)\phi(x)\ \equiv\ \int_\RR e^{ikx} \psi\left(|k|\le k_0\right)\hat\phi(k)\ dk, \label{phi-low}\\
\phi_{high}(x)\ &\equiv\ \left(\ 1\ -\ \psi\left(|D|\le k_0\right) \right)\ \phi(x)\ \equiv\ \int_\RR e^{ikx} \left(\ 1-\psi\left(|k|\le k_0\right)\ \right)\hat\phi(k)\ dk \ .
\label{phi-high}
\end{align}

\subsection{Bounds on $W_+\phi_{low}$}

For $x\ge0$, we can express $W_+\phi_{low}(x)$, in terms of $m_1(x,k)$ which satisfies the bounds \eqref{m1bound}, and for $x\le0$, we can express $W_+\phi_{low}(x)$, in terms of $m_2(x,k)$ which satisfies the bounds \eqref{m2bound}. 
 Since the cases $x\ge0$ and $x\le0$ are very similar, we only carry this calculation out in detail for $x\ge0$. 
   We have, using the notation $P f (x) = f(-x)$, 
\begin{eqnarray*}
&& W_+\phi_{low}\ =\ F^*_{+} F\ \psi(|D|\le k_0) \phi\\
 & = & \int_0^\infty e^{ikx}\ T(k)\  m_1 (x,k) \  \psi\left(|k|\le k_0\right)\ \hat{\phi} (k)\ 
dk +  \int_{-\infty}^0  e^{ikx}\ T(-k)\  m_2 (x,-k) \  \psi\left(|k|\le k_0\right)\  \hat{\phi} (k)\ dk \\
& = & \int_0^\infty e^{ikx}\ T(k)\ m_1 (x,k) \  \psi\left(|k|\le k_0\right)\   \hat{\phi} (k)\ 
dk \\
&& +  \int_{-\infty}^0  e^{ikx}\ [R_1 (-k) e^{-2ikx}\ m_1 (x,-k) + m_1
(x,k) ] \  \psi\left(|k|\le k_0\right)  \hat{\phi} (k)\ dk \\
& = &  \int_0^\infty e^{ikx}\ m_1 (x,k)\ [T(k) + R_1 (k) P] \  \psi\left(|k|\le k_0\right)\  \hat{\phi} (k)\  
dk  +  \int_{-\infty}^0  e^{ikx}\ m_1(x,k)\ \hat{\phi} (k)\ dk, \ \ \ x\ge0,
\end{eqnarray*}
where we have applied \eqref{F+def} and \eqref{Psi+T}.

We continue by using that $\int_0^\infty\ \left[\dots\right]\ dk = \frac{1}{2}\int_{-\infty}^\infty\ (1 + {\rm sgn}(k))\ \left[\dots\right]\ dk$, we have
\begin{eqnarray}
W_+\phi_{low}\ 
& = &\frac{1}{2}\int_{-\infty}^\infty (1 + {\rm sgn}(k))\ e^{ikx} (m_1
(x,k) - 1) T(k) \  \psi\left(|k|\le k_0\right)  \hat{\phi} (k)
dk \label{W+phi}\\
&&+ \frac{1}{2}\int_{-\infty}^\infty (1 + {\rm sgn}(k))\  e^{ikx}  (m_1
(x,k) - 1) R_1 (k) P  \  \psi\left(|k|\le k_0\right) \hat{\phi} (k)
dk \nn\\
&&+\frac{1}{2}\int_{-\infty}^\infty (1 - {\rm sgn}(k))\ e^{ikx}
 (m_1(x,k) - 1 ) \  \psi\left(|k|\le k_0\right)  \hat{\phi} (k) dk \nn \\
&& + \frac{1}{2}\int_{-\infty}^\infty (1 + {\rm sgn}(k))\ e^{ikx}
T(k)  \  \psi\left(|k|\le k_0\right)  \hat{\phi} (k)
dk \nn\\
&& + \frac{1}{2}\int_{-\infty}^\infty (1 + {\rm sgn}(k))\ e^{ikx}  R_1 (k) P \  \psi\left(|k|\le k_0\right) \hat{\phi} (k)
dk \nn\\
&& + \frac{1}{2}\int_{-\infty}^\infty (1 - {\rm sgn}(k))\  e^{ikx}
\  \psi\left(|k|\le k_0\right)  \hat{\phi} (k) dk, \ \ \ x\ge0.
\nn\end{eqnarray}
For $x\le0$ an analogous representation holds with $m_1(x,k)$ replaced by $m_2(x,k)$.

We now show that  $W_{+,1,low}$ is a bounded operator on
$W^{k,p}(\RR_+)$.  
  Each term in the first three lines of \eqref{W+phi} is of the form:
\begin{equation}
\phi\ \mapsto\ S_j\ \circ\ \left(I\pm i\ {\mathcal H}\right)\ \circ\  \Psi(D)\ \phi\ ,
\label{first3}\end{equation}
and each term in the last three lines is of the form
\begin{equation}
\phi\ \mapsto\ \left(I\pm i\ {\mathcal H}\right)\ \circ\  \Psi(D)\ \phi\ ,
\label{last3}\end{equation}
where 
\begin{align}
\Psi(D)\ &=\ F^{-1}\ \hat\Psi(k)\ F\ \ \ {\rm and}\nn\\
\hat\Psi(k) \ &=\ T(k)\ \psi(|k|\le k_0)\ {\rm or}\ R_1(k)\ P\ \psi(|k|\le k_0)\ {\rm or}\ \psi(|k|\le k_0), \nn\\
\left(S_j\Phi\right)(x)\ &\equiv\ \int_\RR\ R_j(x,y)\ \Phi(y)\ dy,\label{Sjdef}\\
R_j(x,y)\ &\equiv\ \int_\RR e^{ikx}\left(m_j(x,k)-1\right)e^{-iky}\ dk.\label{Rjdef}
\end{align}

By hypotheses on $T(k)$ and $R(k)$, $\hat\Psi(k)$ is a multiplier on $W^{k,p}(\RR)$  for $1<p<\infty$ \cite{Stein}.
Therefore, the boundedness of the operators in \eqref{first3} and \eqref{last3}  on $W^{k,p}$ for $1<p<\infty$, and therefore the bound on $W_+\phi_{low}$, follows from 

\begin{lem}
\label{lem:Sjbound}
 $S_1$ is bounded on  $W^{1,p}(\RR_+)$ and  $S_2$ is bounded on  $W^{1,p}(\RR_-)$ for $1<p<\infty$.
 \end{lem}

\noindent {\it Proof of  Lemma \ref{lem:Sjbound}:}
 We focus on the bound for $S_1$ on $W^{k,p}(\RR_+)$. The bound for $S_2$ is bounded on  $W^{k,p}(\RR_-)$ is similar. 

Using the representation formula \eqref{m1B1} we have
\begin{align}
R_j(x,y)\ &\equiv\ \int_\RR e^{ik(x-y)} \int_0^\infty e^{2ikz}\ B_1(x,z)\ dk\ dz\ =\ B_1\left(x,\frac{y-x}{2}\right)
 \nn\end{align}
 and thus the operator $S_1$ simplifies to
\begin{eqnarray*}
\left(S_1\Phi\right)(x) & = &\  \int_{x}^\infty B_1\left(x,\frac{y-x}{2}\right) \Phi(y) dy\ =\
  \int_0^\infty B_1\left(x,\frac{\zeta}{2}\right)\ \Phi(\zeta-x)\ d\zeta,\ \  \ x\ge0\ .
\end{eqnarray*}
Since we must estimate $S_1$ on $W^{1,p}$ we also compute
\begin{eqnarray*}
\partial_x \left(S_1\Phi\right)(x) & = &\  \
  \int_0^\infty B_1\left(x,\frac{\zeta}{2}\right)\ (-\partial_\zeta)\Phi(\zeta-x)\ d\zeta\ +\ \ 
 \int_0^\infty \partial_x B_1\left(x,\frac{\zeta}{2}\right)\ \Phi(\zeta-x)\ d\zeta\\
 & = &   \int_x^\infty B_1\left(x,\frac{y-x}{2}\right)\ (-\partial_y)\Phi(y)\ dy\ +\ \ 
 \int_x^\infty \partial_x B_1\left(x,\frac{y-x}{2}\right)\ \Phi(y)\ dy,\ \   \ x\ge0\ .
\end{eqnarray*}

To prove boundedness of $S_1 $ and $\partial S_1$ on  $L^p$ of the operator we use 
that the operator
\begin{equation}
S_R\Phi(x)\ =\ \int_\RR R(x,y)\ \Phi(y)\ dy, 
\nn\end{equation}
is bounded on $L^p$ with estimate
\begin{equation} \left\| S_R\Phi\ \right\|_{L^p}\ \le\ C_R\ \left\| \Phi\right\|_{L^p}, \ \ 1\le p\le\infty
\label{youngineq}
\end{equation}
 if
\begin{equation}
C_R\equiv \sup_{x\ge0}\ \int_\RR\ |R(x,y)|\ dy\ +\ \sup_{y\ge0}\ \int_\RR |R(x,y)|\ dx\ <\ \infty\ .
\label{youngprime} \end{equation}

Note that by \eqref{eqn:b1bd} and \eqref{eqn:b1primebd} we have
\begin{eqnarray}\label{b1-bounds}
\left|B_1 (x,z)\right| & \lesssim & \int_{x+z}^\infty |V_{reg}(s)| ds\ \ {\rm and}\ \ 
\left|\partial_x B_1 (x,z)\right|  \lesssim  |V_{reg}(x)|\ +\  \int_{x+z}^\infty |V_{reg}(s)| ds .
\end{eqnarray}
Therefore,
\begin{eqnarray*}
 && \!\!\!\!\!\!\!\!\!\!\!\!\!  \sup_{x\ge0} \int  1_{y\ge x}\left| B_1\left(x,\frac{y-x}{2}\right)\right|\ dy\ +\ 
\sup_{y\ge0} \int  1_{y\ge x}\left|
  B_1\left(x,\frac{y-x}{2}\right)\right|\ dx\nn \\
&&\le \ 
2\sup_{x\ge0} \int_0^\infty \int_{\frac{x+y}{2}}^\infty\ |V_{reg}(s)| ds\ dy\nn\\
&& \le \ 
2\int_0^\infty\ \left(1+\frac{x+y}{2}\right)^{-\frac{3}{2}-}\int_{\frac{x+y}{2}}^\infty(1+s)^{{3\over2}+}|V_{reg}(s)| ds\ \\
&& \le const\times \|V_{reg}\|_{L^1_{{3\over2}+}(\RR)}.\nn
\end{eqnarray*}
A similar bound applies to the kernel $1_{x\ge y}\partial_x B_1\left(x,\frac{y-x}{2}\right)$. 
 Thus, we have 
 \begin{equation}
\| S_1\Phi\|_{W^{1,p}(\RR_+)}  \equiv \| S_1\Phi\|_{L^p(\RR_+)} + \|\partial_x \left(S_1\Phi\right)\|_{L^p(\RR_+)} \le\ C\|V_{reg}\|_{L^1_{{3\over2}+}(\RR)}\ \|\Phi\|_{W^{1,p}(\RR_+)}\ .\nn\end{equation}
Applying similar arguments with $S_1$ replaced by $S_2$ for $x\le0$ yields boundedness of $S_2$ on $W^{1,p}$, from which we conclude
\begin{equation}
\| W_+\phi_{low}\|_{W^{1,p}(\RR)} \ \le\ C\ \|V_{reg}\|_{L^1_{{3\over2}+}(\RR)}\ \|\phi\|_{W^{1,p}(\RR)}\ .
\label{low-bound}
\end{equation}
 This completes the low frequency analysis.

\subsection{High Frequencies}

We have, using \eqref{eqn:f1f2} and the notation $ P f (x) = f(-x)$, 
\begin{eqnarray*}
W_+\phi_{high} &=& F^*_{+} F \left(1-\psi(|D|\le k_0)\right) \phi\\
 & = & \int_0^\infty T(k) e^{ikx} m_1 (x,k) (1- \psi\left(|k|\le k_0\right)) \hat{\phi} (k)
dk \\
&& +  \int_{-\infty}^0  T(-k) e^{ikx} m_2 (x,-k) (1-\psi\left(|k|\le k_0\right))  \hat{\phi} (k) dk \\
& = & \int_0^\infty T(k) e^{ikx} m_1 (x,k) (1-\psi\left(|k|\le k_0\right))   \hat{\phi} (k)
dk \\
& + & \int_{-\infty}^0  e^{ikx} [R_1 (-k) e^{-2ikx} m_1 (x,-k) + m_1
(x,k) ] (1-\psi\left(|k|\le k_0\right))  \hat{\phi} (k) dk \\
& = &  \int_0^\infty e^{ikx} m_1 (x,k) [T(k) + R_1 (k) P] (1-\psi\left(|k|\le k_0\right))  \hat{\phi} (k)
dk   +  \int_{-\infty}^0  e^{ikx} m_1(x,k) \hat{\phi} (k) dk.
\end{eqnarray*}

For $x\ge0$ we rewrite this expression as 
\begin{eqnarray*}
W_+\phi_{high} & = &\frac{1}{2}\int_{-\infty}^\infty e^{ikx}\ (1 + {\rm sgn}(k)) \ (m_1
(x,k) - 1) T(k) (1-\psi\left(|k|\le k_0\right))  \hat{\phi} (k)
dk \\
&&+ \frac{1}{2}\int_{-\infty}^\infty  e^{ikx}\ (1 + {\rm sgn}(k))  (m_1
(x,k) - 1) R_1 (k) P  (1-\psi\left(|k|\le k_0\right)) \hat{\phi} (k)
dk \\
&&+  \frac{1}{2}\int_{-\infty}^\infty  e^{ikx} \ (1 - {\rm sgn}(k))  
 (m_1(x,k) -1 ) (1-\psi\left(|k|\le k_0\right))  \hat{\phi} (k) dk \\
&& +  \frac{1}{2}\int_{-\infty}^\infty e^{ikx}\ (1 + {\rm sgn}(k)) 
 T(k)  (1-\psi\left(|k|\le k_0\right))  \hat{\phi} (k)
dk \\
&&+  \frac{1}{2}\int_{-\infty}^\infty e^{ikx}  \ (1 + {\rm sgn}(k))   R_1 (k) P (1-\psi\left(|k|\le k_0\right)) \hat{\phi} (k)
dk \\
&& +  \frac{1}{2}\int_{-\infty}^\infty  e^{ikx}\ (1 - {\rm sgn}(k)) 
 (1-\psi\left(|k|\le k_0\right))  \hat{\phi} (k) dk ,\ \ x\ge0.
\end{eqnarray*}
An analogous expression, with $m_1(x,k)$ replaced by $m_2(x,k)$, is used  for $x\le0$.
 We proceed now to show that each term is bounded on $W^{1,p}(\RR_+)$, $p\ge1$.
 
 Each summand in this decomposition of $W_+\phi_{high}$ is of the form:
 \begin{align}
 &\phi\mapsto S_j\circ \rho(D)\ \phi,\ \ \ \ \ \ {\rm or}\ \ \ \ \phi\mapsto \rho(D) \phi.
 \label{model-high}
 \end{align}
 where $\rho(D)=F^{-1} \hat\rho(k) F$. 
Here, $S_j,\ j=1,2$,  defined  in \eqref{Sjdef} and \eqref{Rjdef},  is bounded on $W^{1,p}(\RR_+)$ for $1<p<\infty$, as proved in the previous section. Moreover, 
 $\rho(k)$ is a multiplier on $W^{1,p}(\RR)$ for $1<p<\infty$ due to hypotheses on $R(k), T(k)-1, \partial_k R(k)$ and $\partial_k T(k)$, and the fact that $1-\psi(|k|\le k_0)$ is smooth, asymptotically constant as $k\to\infty$ and vanishing in a neighborhood of $0$.
 It follows that 
 \begin{equation}
 \| W_+\phi_{high}\|_{W^{1,p}(\RR_+)}\le C \|V_{reg}\|_{L^{1}_{{3\over2}+}(\RR)}\ \|\phi\|_{W^{1,p}(\RR_+)}.\label{W+highbound}
 \end{equation}
 An estimate analogous to \eqref{W+highbound}, similarly proved using a representation of $W_+\phi_{high}(x)$ for $x\le0$, in terms of $S_2$, also holds. Thus,
 \begin{equation}
 \| W_+\phi_{high}\|_{W^{1,p}(\RR)}\le C \|V_{reg}\|_{L^{1}_{{3\over2}+}(\RR)}\ \|\phi\|_{W^{1,p}(\RR)}\ .
\label{high-bound}
\end{equation}

The decomposition \eqref{low+high} and the bounds \eqref{low-bound} and \eqref{high-bound} imply the result. This completes the proof of the main result,  Theorem \ref{thm:main}.

\section{Examples and Applications}
\label{sec:nonapp}

\subsection{$V(x)=$ a sum of Dirac delta masses}
\label{sec:singpot}

In this section we verify hypotheses \eqref{RT-assume},
\eqref{m1bound}, \eqref{m2bound} for the case
 of a potential, which is the sum of Dirac delta functions, thereby establishing the applicability of our main results to this case.

 We follow the analysis from \cite{HMZ} and \cite{TZ}, see also \cite{GS},
  \cite{GT} for specific examples.   Seek solutions of the form
\begin{eqnarray}
\left(\ H_{\vec{q},\vec{y}} -  {1\over2}k^2\ \right) e_\pm (x,k) = 0,
\end{eqnarray}
where $H_{\vec{q},\vec{y}}=\sum_{j=0}^{N-1} q_j \delta(x-y_j)$ when $\vec{q}=(q_0,\cdots,q_{N-1}),\ \vec{y}=(y_0,\cdots,y_{N-1})$, and where $e_{\pm}(x,k)$ represent the  distorted Fourier basis functions as defined \eqref{outgoing}.  Thus,\begin{eqnarray}
e_+ (x,k) = \left\{ \begin{array}{c}
e^{i k x} + B_0 e^{-i k x} \ \text{for $x < y_0$},  \\
A_1 e^{i k x} + B_1 e^{-i k x} \ \text{for $y_0 < x < y_1$}, \\
\vdots \\
A_N e^{i k x}  \ \text{for $x > y_{N-1}$}, 
\end{array} \right. 
\end{eqnarray}
where we have taken $A_0 = 1$ and $B_N = 0$.  With this choice of notation, we have, referring to \eqref{u1exp} and \eqref{u2exp}, $A_N = T$ the
transmission coefficient and $B_0 = R_1$ the reflection coefficient
for the ``incoming'' plane wave $e^{i k x}$ from $-\infty$.  Then, we have the following system of equations implied by continuity and jump conditions at the points $\{y_j\}$ for $j = 0, \dots, N-1$:
\begin{eqnarray*}
e^{i k y_0} + B_0 e^{- i k y_0} & = & A_1 e^{i k x_0} + B_1 e^{- i k y_0} \\
i k \left[ A_1 e^{i k y_0} - B_1 e^{- i k y_0} -
  e^{i k y_0} + B_0 e^{- i k y_0} \right] & = & 2 q_0
\left[ A_1 e^{i k y_0} + B_1 e^{- i k y_0} \right]\\
\vdots \\
A_{N-1} e^{i k y_{N-1}} + B_{N-1} e^{- i k y_{N-1}} & = & A_N e^{i k y_{N-1}} \\
i k \left[ A_{N} e^{i k y_{N-1}} - A_{N-1} e^{i k
    y_0} + B_{N-1} e^{- i k y_0} \right] & = & 2 q_{N-1} \left[
  A_N e^{i k y_{N-1}} \right] .
\end{eqnarray*}
Note, the above system guarantees unitarity, or that
\begin{eqnarray}
|B_0|^2 + |A_N|^2 = 1.
\end{eqnarray}

We can define similarly
\begin{eqnarray}
e_- (x,k) = \left\{ \begin{array}{c}
D_0 e^{-i k x} \ \text{for $x < y_0$},  \\
C_1 e^{i k x} + D_1 e^{-i k x} \ \text{for $y_0 < x < y_1$}, \\
\vdots \\
C_N e^{i k x} +  e^{-i k x} \ \text{for $x > y_{N-1}$},
\end{array} \right. 
\end{eqnarray}
where now the incoming wave is $e^{-i k x}$ from $\infty$ and
the scattering matrix is determined by the transmission coefficients
$D_0 = T$ and the reflection coefficient $C_N = R_2$
for the ``incoming'' plane wave $e^{-i k x}$ from $\infty$.  

\subsubsection{Bounds on $m_1$, $m_2$:}

In addition, for general singular potentials with compact support, we have 
\begin{eqnarray*}
m_1 (x,k) \!\!\!\!\! & = & \!\!\!\!\!  e^{-ikx} f_1 (x,k) =
\left\{ \begin{array}{l}
e^{-ikx} \frac{e_{+} (x,k)}{T(k)} \ \text{for} \ x < y_{N-1} , \\
1,\ \text{for} \ x > y_{N-1} ,
\end{array} \right. \\
m_2 (x,k) \!\!\!\!\! & = & \!\!\!\!\!   e^{ikx} f_2 (x,k) =\left\{ \begin{array}{l}
e^{ikx} \frac{e_{-} (x,k)}{T(k)} \ \text{for} \ x > y_{0} , \\
1,\ \text{for} \ x < y_{0} .
\end{array} \right.
\end{eqnarray*}
Hence, there exists constants $C_\alpha^1 (y_{N-1})$ and $C_\alpha^2 (y_{0})$ such that
\begin{align}
\label{eqn:regm1}
|\partial_k^\alpha m_1 (x,k)  | &\leq C_\alpha^1 (y_{N-1}) \ \text{for}
\ y_{N-1}  > x \geq 0, \\
\label{eqn:regm2}
|\partial_k^\alpha m_2 (x,k)  | &\leq  C_\alpha^2 (y_{0}) \ \text{for}
\ y_{0} < x \leq 0.
\end{align}
As a result, we see that an arbitrary collection of $\delta$ functions
satisfies assumptions \eqref{m1m2compbound} , \eqref{m1bound} and \eqref{m2bound} as
required for the proof of Theorem \ref{thm:main}.

We conclude this subsection with explicit computations of the transmission and reflection coefficients for single and double $\delta$ well potentials:

\subsubsection{{\bf Single $\delta$ potential ($H_{q} = - q \delta(x) $)}:}

Setting up the appropriate equations, we have
\begin{eqnarray}
\label{eqn:refco-d}
R_1 = r_{q} \!\!\!\! & = & \!\!\!\!\! \frac{q}{i k -q}, \\
\label{eqn:tranco-d}
T = t_{q} \!\!\!\! & = & \!\!\!\!\! \frac{i k}{i k -q},
\end{eqnarray}
where $r_{q}$, $t_{q}$ are the reflection and
transmission coefficients for $H_q$ respectively.  We must show the bounds
from \eqref{RT-assume} hold, however such bounds follow clearly for \eqref{eqn:tranco-d}, \eqref{eqn:refco-d}.

\subsubsection{{\bf Double $\delta$ potential ($H_{q,L} = - q ( \delta(x+L) + \delta(x-L) )$)}:}

Setting up the appropriate equations, we have
\begin{eqnarray}
\label{eqn:refco-dd}
R_1 = r_{q,L} & = & \left( \frac{q(ik -q) e^{2ik L} + q (ik+q) e^{-2ik L}}{q^2 e^{2i k L} - (i k +q)^2 e^{-2ik L}} \right) e^{-2i k L}, \\
\label{eqn:tranco-dd}
T = t_{q,L} & = & \left( \frac{k^2}{q^2 e^{2i k L} - (i k +q)^2 e^{-2ik L}} \right) e^{-2i k L},
\end{eqnarray}
where $r_{q,L}$, $t_{q,L}$ are the reflection and
transmission coefficients for $H_{q,L}$ respectively.  

Again, we must verify bounds \eqref{RT-assume}, hence we must prove
for instance
\begin{eqnarray*}
\dot{t}_{q,L} (k) \leq C (1 + |k|)^{-1}  ,
\end{eqnarray*}
provided $qL \neq 1/2$.  Indeed,
we have
\begin{eqnarray*}
\dot{t}_{q,L} (k) & = & \frac{2 k (k^2 - 2i k q + q^2
  (e^{4i k L} -1)) - 2i k^2 ( 2 L q^2 e^{4i k L} -
  (i k + q))}{(k^2 - 2i k q + q^2
  (e^{4i k L} -1))^2},
\end{eqnarray*}
which satisfies
\begin{eqnarray*}
|\dot{t}_{q,L} (k)| \sim  \mathcal{O} (|k|^{-1})
\end{eqnarray*}
as $k \to \infty$ and 
\begin{eqnarray*}
|\dot{t}_{q,L} (k)| \sim  \mathcal{O} (\frac{1}{4 q^2 L - 2q})
\end{eqnarray*}
as $k \to 0$.  A similar computation holds for $r_{q,L}$.

\begin{rem}
Such bounds can be verified for a more general $\delta$ function
potential using the expressions 
\begin{eqnarray*}
T(k) & = & 1 + \frac{\int_{-\infty}^\infty V(t) dt}{2ik} + \mathcal{O}
(k^{-2}), \\
R_j(k) & = & \frac{T(k) \int_{-\infty}^\infty e^{\pm 2ikt} V(t) dt}{2ik} + \mathcal{O}
(k^{-2}),
\end{eqnarray*}
which can be derived from the expressions for $m_1,m_2$ as in
\cite{DT}.
\end{rem}

\subsection{Commutator / Resolvent type bounds}
\label{sec:res}

In \cite{DucheneWeinstein}, where homogenization of high contrast oscillatory structures with defects is studied, bounds on $(H_0  + 1)^{-1} (H_{\vec{q},\vec{y}} + 1)$ are required to estimate a Lipmann Schwinger equation.  We have, by our main theorem that 
\begin{eqnarray*}
(H_0  + 1)^{-1} (H_{\vec{q},\vec{y}} + 1) P_c = (H_0 +1)^{-1} W_+ (H_0 +1) W_+^* : L^2 \to L^2.
\end{eqnarray*}

\subsection{Dispersive and Strichartz estimates in $H^1$ for $\delta$-NLS}
\label{strich}

We may represent 
\begin{eqnarray}
e^{-i t H} P_c f =  \frac{1}{2\pi} \int \int_0^\infty e^{-\frac{it k^2}{2}} \left( e_+(x,k)\ \overline{e_+ (x,k)} + 
e_- (x,k)\ \overline{e_- (x,k)} \right) f(y) dk dy.
\end{eqnarray}
From here, we may use direct computations to arrive at Strichartz estimates and apply Weder's results on wave operators since the potentials are all in $L^1$ with compact support. 

Using the properties of wave operators, we have
\begin{eqnarray}
\| e^{i H t} P_c f \|_{L^p} = \| W_\pm e^{i tH_0} W^*_\pm f \|_{L^p}
\end{eqnarray}
and using standard dispersive estimates for the linear Schr\"odinger operator (see for instance \cite{SS} for a concise overview) arrive at
\begin{eqnarray}
\| e^{i H t} P_c f \|_{L^p} \leq C_p t^{-(\frac{1}{2} - \frac{1}{p})} \| f \|_{W^{1,p}}.
\end{eqnarray}
Define a Strichartz pair $(q,r)$ to be admissible if
\begin{eqnarray}
\frac{2}{q} = \frac{1}{2} - \frac{1}{r}
\end{eqnarray}
with $2 \leq r < \infty$.  Then, we arrive at the celebrated Strichartz estimates
\begin{eqnarray}
\label{eqn:strich1}
\| e^{iHt} P_c u_0 \|_{L^q W^{1,r}} \lesssim \| u_0 \|_{ W^{1,2}}
\end{eqnarray}
and
\begin{eqnarray}
\label{eqn:strich2}
\left\| \int_0^t e^{iH(t-s)} P_c f \right\|_{L^q W^{1,r}} \lesssim \| f(x,t) \|_{L^{\tilde{q}}_t W^{1,\tilde{r}}_x}
\end{eqnarray}
using duality techniques and once again the boundedness of the wave operators.

As a side note, using positive commutators and well crafted local smoothing spaces, from \cite{MMT} we have the Strichartz estimate
\begin{eqnarray}
\label{eqn:str1}
\left\| \int_0^t e^{iH(t-s)} P_c f  \right\|_{L^\infty L^2} \lesssim \| f(x,t) \|_{L^{\tilde{p}}_t L^{\tilde{q}}_x}.
\end{eqnarray}
Now, by boundedness of wave operators on $W^{1,p}$ spaces for singular potentials as proved in Theorem \ref{thm:main}, we have the following useful relation
\begin{eqnarray}
\left\| \int_0^t e^{iH(t-s)} P_c f  \right\|_{L^\infty H^1} \lesssim \| f(x,t) \|_{L^{\tilde{p}}_t W^{1,\tilde{q}}_x},
\end{eqnarray}
where $(\tilde{p},\tilde{q})$ is a dual Strichartz pair without first going through the dispersive estimates.

\subsection{Local Well-Posedness in $H^1$ for $\delta$-NLS}
\label{sec:lwp}

Consider the nonlinear Schr\"odinger / Gross-Pitaevskii, with a potential consisting of a finite 
set of Dirac delta functions:
\begin{eqnarray*}
\left\{ \begin{array}{c}
i \partial_t u + H_{\vec{q},\vec{y}} u - |u|^{2\sigma} u = 0 , \\
u(x,0) = u_0 (x) \in H^1,
\end{array} \right.
\end{eqnarray*}
for $0 < \sigma < \infty$.  We seek a solution in the following sense:
\begin{equation}
u=\Lambda[u],
\nn\end{equation} where 
\begin{eqnarray}
\Lambda[u](t) = e^{-i H_{\vec{q},\vec{y}} t} u_0  - i \int_0^t e^{-i
  H_{\vec{q},\vec{y}}( t-s)} |u|^{2\sigma} u (s) ds.
  \label{nls-int}
\end{eqnarray}
We claim that local well-posedness  can be established via the contraction mapping principle in the space
 $C^0([0,T);H^1(\RR))$ for $T$ sufficiently small. 
  To prove the necessary boundedness and contraction estimates, it is natural to 
  apply the operator $(I+H_{\vec{q},\vec{y}})^{1\over2}P_c$, which commutes with 
   the group $e^{-i H_{\vec{q},\vec{y}} t} $ to \eqref{nls-int}. Then, estimates follow in a straightforward way, using  that $H^1(\RR)$ is an algebra, provided the space 
   \begin{equation}
 \mathcal{H}^1(\RR)=\left\{ f: (I+H_{\vec{q},\vec{y}})^{1\over2}P_c f\in L^2(\RR)\ \right\}
 \end{equation}
 is equivalent to the classical Sobelev space $H^1$.  This follows from the relations
 \begin{equation}
 (I+H)^{1\over2}P_c = W(I-\partial_x^2)^{1\over2}W^*,\ \ W^*(I+H)^{1\over2}W = (I-\partial_x^2)^{1\over2}
 \nn\end{equation}
 and our results on the boundedness of wave operators associated with $H_{\vec{q},\vec{y}}$ on $H^1$.

\subsection{Long time dynamics for $NLS$ with a double $\delta$ well
  potential}
\label{sec:dwp}

In \cite{MW}, the long time dynamics of solutions to the
nonlinear Schr\"odinger / Gross-Pitaevskii equation
\begin{eqnarray}
\label{eqn:nlsdwp-gK}
i \partial_t u =  (- \Delta + V(x)) u\ +\ gK\left[\left|u\right|^2\right] u,
\end{eqnarray}
where $V$ is a symmetric, double well potential, are studied.
In particular, under appropriate spectral assumptions on the operator $H
=- \partial_x^2 + V(x)$, in a neighborhood of a symmetry breaking bifurcation point, 
 there are different classes of oscillating solutions \eqref{eqn:nlsdwp-gK} which {\it shadow} periodic orbits of a finite dimensional reduction on very long, but finite, time scales. These solutions correspond to states 
with mass concentrations oscillating between the two wells of a symmetric potential well.   The proof  requires dispersive / Strichartz type estimates. 
The results of this paper imply that the results of \cite{MW} extend to \eqref{eqn:nlsdwp-gK}
for  the case of singular potentials, such as 
\begin{eqnarray*}
V(x) = -q[ \delta(x-L) + \delta (x+L)].
\end{eqnarray*}

\appendix

\section{Bounds on  $m_j(x;k),\ j=1,2$}\label{app:mjbounds}
 Denote by 
$m_1(x,k) = e^{-ikx} f_1 (x,k)$ and $m_2 (x,k) = e^{ikx} f_2 (x,k)$.  Then, we have
\begin{eqnarray*}
m_1 (x,k) & = & 1 + \int_x^\infty D_k (y-x) V(y) m_1 (y,k) dy, \\
m_2 (x,k) & = & 1 + \int_{-\infty}^x D_k (x-y) V(y) m_2 (y,k) dy, \\
D_k (x) & = & \int_0^x e^{2iky} dy.
\end{eqnarray*}

We remark the derivation for $m_1(x,k),\ x\ge0$. Similar remarks apply to $m_2(x,k)$ on $x\le0$. 
By results in  \cite{DT}, for $V\in L^1_{{3\over2}+}(\RR)$ the function  $m_1 (x,k) -1$ is in the Hardy space, and therefore 
there exists $B_1 \in L^2 (\RR_+)$ such that
\begin{eqnarray}
m_1 (x,k) = 1 + \int_0^\infty B_1 (x,y) e^{2iky} dy.\ \label{m1B1}
\end{eqnarray}
Moreover, 
\begin{align}
\label{eqn:b1bd}
| B_1 (x,y)| &\le\ C\ e^{\gamma_1 (x)}\ \int_{x+y}^\infty |V(t)| dt, \ x, \
  y > 0, \\
\gamma_1 (x) &= \int_x^\infty (t-x) |V(t)| dt.
\end{align}
Similarly, 
\begin{align}
\label{eqn:b1primebd}
\left| \partial_x B_1 (x,y)\right| &\leq \ C\ e^{\gamma_1 (x)}\ 
\left(\  V(x+y)\ +\ 
\int_{x+y}^\infty |V(t)| dt\ \right) , \ x \in \mathbb{R}, \
  y > 0, \\
\gamma_1 (x) &= \int_x^\infty (t-x) |V(t)| dt.
\end{align}

The proof of \cite{DT} extends to the case where $V(x)=V_{sing}(x)+V_{reg}(x)$, where $V\in L^1_{{3\over2}+}(\RR)$  and $V_{sing}(x)$ consists of a finite sum of delta functions. Indeed,  for $V (x) = \delta (x)$ we have
\begin{eqnarray*}
B_1 (x,y) = \sum_{n=0}^\infty K_n (x,y)
\end{eqnarray*}
for 
\begin{eqnarray*}
K_0 (x,y) = \int_{x+y}^\infty V(t) dt, \ K_{n+1} (x,y) = \int_0^y
\int_{x+y-z}^\infty V(t) K_n (t,z) dt dz, \ n = 0,1,\dots .
\end{eqnarray*}
Hence,
\begin{eqnarray*}
K_0  = \int_{x+y}^\infty \delta (t) dt = \left\{ \begin{array}{c}
1, \ x+y < 0 \\
\frac{1}{2}, \ x+y = 0 \\
0, \ x+y > 0.
\end{array} \right.
\end{eqnarray*}
As a result,
\begin{eqnarray*}
K_1 (x,y) & = & \int_0^y \left\{ \begin{array}{c}
0, \ x+y-z > 0 \\
\frac12 K(0,z), \ x+y -z = 0 \\
K(0,z), \ x+y-z < 0
\end{array} \right. \\
& = & 0
\end{eqnarray*}
since $K(0,z) = 0$ for any $z > 0$.  Similar computations can be done for
larger collections of $\delta$ functions.
Hence, for $V = \delta$, we have
\begin{eqnarray*}
B_1 (x,y) = K_0 (x,y)
\end{eqnarray*}
for which the bounds \eqref{eqn:b1bd}, \eqref{eqn:b1primebd} hold
obviously in the sense of distributions.  \


\begin{thebibliography}{MA}

\bibitem{Ag} S. Agmon.  {\em Spectral properties for Schr\"odinger operators and scattering theory,} {\it Ann. Scuola Norm. Sup. Pisa Cl. Sci. (4)}, {\bf 2} (1975), 151-218.

\bibitem{Caz} T. Cazenave.  {\em Semilinear Schr\"odinger Equations},
  Courant Lecture Notes in Mathematics, Vol. {\bf 10}.  New York
  University, Courant Institute of Mathematical Sciences, New York;
  American Mathematical Society, Providence, RI (2003).

\bibitem{DF} P. D'Ancona and L. Fanelli.  {\em $L^p$ boundedness of the
    wave operators for the one dimensional Schr\"odinger operator,} {\it Communications in Mathematical Physics}, {\bf 268}, Number 2 (2006), 415-438.

\bibitem{DT} P. Deift and E. Trubowitz.  {\em Inverse scattering on
    the line}, {\it Commun. Pure Appl. Math.}, {\bf 32} (1979), 121-251.
    
  \bibitem{DucheneWeinstein} V. Duch\^ene and M.I. Weinstein.  {\em Scattering and homogenization for truncated periodic structures with defects}, in preparation.
    
    \bibitem{Folland} G.B. Folland  {\it  Partial Differential Operators }, Princeton University Press 


\bibitem{GHW} R. Goodman, P. Holmes and M.I. Weinstein.  {\em Strong NLS
    soliton-defect interactions}, {\it P hys. D}, {\bf 192}, No. 3-4
  (2004), 215-248.
  
  \bibitem{GS} D.J. Griffiths and C.A. Steinke.  {\em Waves in locally periodic media}, {\it Am. J. Physics}, {\bf 69 }, No. 2
  (2000), 137--154.

  \bibitem{GT} D.J. Griffiths and N.F. Taussig.  {\em Scattering from
      a locally periodic potential}, {\it Am. J. Physics}, {\bf 60 }, No. 10
  (1992), 883--888.

\bibitem{HMZ} J. Holmer, J.L. Marzuola and M. Zworski.  {\em Fast soliton scattering by delta impurities,} {\it Communications in Mathematical Physics}, {\bf 274}, Number 1 (2007), 187-216.

\bibitem{HMZ1} J. Holmer, J.L. Marzuola and M. Zworski.  {\em Soliton splitting by external delta potentials,} {\it Journal of Nonlinear Science}, {\bf 17}, Number 4 (2007), 349-367.

\bibitem{HZ}  J. Holmer and M. Zworski.  {\em Slow soliton interaction
    with delta impurities}, {\it J. Mod. Dyn. 1}, No. 4 (2007), 689-718.

\bibitem{Ho2} L. H\"ormander.  {\it The Analysis of Linear Partial Differential Operators II}, Classics in Mathematics.  Springer-Verlag, Berlin (2005).

\bibitem{JW} R.K. Jackson and M.I. Weinstein.  {\em Geometric analysis of bifurcation and symmetry breaking in a Gross-Pitaevskii Equation}, {\it Journal of Statistical Physics}, {\bf 116}, No. 1 (2004), 881-905.

\bibitem{KKSW} P. Kevrekidis, E. Kirr,  E. Shlizerman and M.I. Weinstein. {\em Symmetry breaking bifurcation in nonlinear Schr\"odinger/Gross-Pitaevskii equations}, {\it SIAM J. Math. Anal.}, {\bf 40}, No. 2 (2008), 566-604.  

\bibitem{KS} J. Krieger and W. Schlag.  {\em Stable manifolds for all
    monic supercritical focusing nonlinear Schr\"odinger equations in
    one dimension}, {\it J. Amer. Math. Soc.}, {\bf 19}, No. 4 (2006),
  815-920.

\bibitem{MW} J.L. Marzuola and M.I. Weinstein.  {\em Long time
    dynamics near the symmetry breaking bifurcation for Nonlinear
    Schr\"odinger/Gross-Pitaevskii Equations}, to appear in {\it
    Disc. and Cont. Dyn. Syst. A}.

\bibitem{MMT} J.L.  Marzuola , J. Metcalfe, and D. Tataru. {\em Strichartz
    estimates and local smoothing estimates for asymptotically flat
    Schr\"odinger equations}, {\it J. Funct. Anal.}, {\bf 255}, Issue
  6 (2008), 1497-1553.

\bibitem{Sad} C. Sadowsky.  {\it Interpolation of Operators and
    Singular Integrals}. New York:  Marcel Dekker (1979).

\bibitem{RSv3} M. Reed and B. Simon. {\it Methods of modern
    mathematical physics. III. Fourier Analysis}, Academic Press, New York-London (1978).

\bibitem{RSv4} M. Reed and B. Simon. {\it Methods of modern mathematical physics. IV. Analysis of Operators}, Academic Press, New York-London (1978).

\bibitem{Schec} M. Schecter. {\it Operator Methods in Quantum
    Mechanics}, New York:  North Holland (1979).

\bibitem{Sch} W. Schlag. {\em Spectral theory and nonlinear partial differential equations: a survey}, {\it Disc. and Cont. Dyn. Syst.}, {\bf 15}, No. 3 (2006), 703-723.

\bibitem{Stein} E.M. Stein.  {\it Singular Integrals and Differentiability Properties of Functions}, Princeton University Press (1970)

\bibitem{SS} C. Sulem and P. Sulem. {\it The Nonlinear Schrodinger Equation. Self-focusing and wave-collapse}, Applied Mathematical Sciences, {\bf 39}. Springer-Verlag, New York (1999).

\bibitem{TZ} S. H. Tang and M. Zworski. {\it Potential Scattering on the Real Line}, unpublished lecture notes.

%\bibitem{VNMA}  {\bf XYZ???} P. Varatharajah, A. Newell, J. Moloney and A. Aceves.  {\em Transmission, reflection and trapping of collimated light beams in diffusive Kerr-like nonlinear media}, {\it Phys. Rev. A.}, {\bf 42} (1990), 1767-1774.

\bibitem{W} R. Weder.  {\em The $W_{k,p}$-Continuity of the
    Schr\"odinger Wave Operators on the Line}, {\it
    Commun. Math. Phys.}, {\bf 208} (1999), 507-520.

\bibitem{Ya} K. Yajima.  {\em The $W^{k,p}$-continuity of wave
    operators for Schr\"odinger operators}, {\it J. Math. Soc. Japan},
  {\bf 47} (1995), 551-581.

\end{thebibliography}
\end{document}